\documentclass[12pt]{amsart}
\usepackage{amsmath,verbatim}
\numberwithin{equation}{section}
\def\setm{\begin{picture}(12,12)\put(1,6){\line(5,-2){10}}\end{picture}}

\newtheorem{thm}{Theorem}[section]
\newtheorem{lem}[thm]{Lemma}
\newtheorem{cor}[thm]{Corollary}
\newtheorem{prop}[thm]{Proposition}

\theoremstyle{definition}
\newtheorem{defn}{Definition}[section]   
\newtheorem{conj}{Conjecture}[section]	 
\newtheorem{exmp}{Example}[section]

\theoremstyle{remark}
        
\newtheorem{note}{Note}         

\newcommand{\R}{{\mathbb R}} 
\newcommand{\Z}{{\mathbb Z}}
\def\approaches{\rightarrow} \def\cross{\times} \def\ident{\equiv}
\def\takes{\colon} 
\let\limsup\varlimsup

\date{November 14, 1996; revised January 17, 1997}
\thanks{The first author was funded in part by NSF grant DMS 94-04278,
and both authors enjoyed the hospitality of the IMA during the
Spring 1996 term}

\title{On Distortion and Thickness of Knots}

\author{Robert B. Kusner}
\address{Institute for Advanced Study}
\email{kusner@math.umass.edu}

\author{John M. Sullivan}
\address{University of Illinois}
\email{sullivan@geom.umn.edu}

\begin{document}
\maketitle

\section*{Introduction}

What length of rope (of given diameter) is required to tie a particular knot?
Or, to turn the problem around, given an embedded curve,
how thick a regular neighborhood of the curve also is embedded?
Intuitively, the diameter of the possible rope is bounded by
the distance between strands at the closest crossing in the knot. 
But of course the distance between two points along a curve goes
to zero as the points approach each other, so to make the notion
precise, we need to exclude some neighborhood of the diagonal.

Various notions of thickness have been proposed recently.  For
example, \cite{Lith1} defines a thickness by considering the distance
function between points on the curve only where it has critical
points.  But the definition there also involves the minimum radius of
curvature of the curve, and thus is unbounded for polygonal curves.
Here we introduce and compare two new families of thickness measures.
One makes use of Gromov's concept of distortion (see \cite{GromovHED}
and \cite[p.~114]{GromovFRM} and \cite[pp.~6--9]{GLP}); it applies to
all rectifiable curves (including polygons).  The other generalizes
the notion from \cite{Lith1}.  Our main result is a basic inequality
(Theorem~\ref{distocurv}) between these measures.

The distortion thickness should permit us to prove the existence of thickest
curves of prescribed length (or dually, shortest curves of prescribed
thickness) in each knot class; such curves are of interest to
chemists and biologists modeling polymers and DNA (see, for example,
\cite{Katritch}).  Moreover, curvature bounds should follow from the
optimality, and we conjecture that, for an optimal knot, all
reasonable measures of thickness should be equal.

\section{Notation and Definitions}
We will deal throughout with an embedded rectifiable closed curve
$\gamma$ in $\R^n$, of finite length $L(\gamma)$.
Thus we may assume that $\gamma\takes \R/L\Z \to \R^n$ is a
Lipschitz parameterization by arclength.
(Often we will normalize by rescaling so that $L=2\pi$.)
If $p,q \in \gamma$ are two points on the curve, $|p-q|$ will denote
the straight-line (chord) distance between them in $\R^n$, and $d(p,q)$
the arclength distance between them along $\gamma$.  This arclength
distance is always measured the shorter way around $\gamma$, so that
it is never more than $L/2$.  Given any point $p$ on $\gamma$, there
is a unique {\em opposite point} $p^*$ such that $d(p,p^*)=L/2$.

\begin{defn}
The {\em distortion} between distinct points $p$ and $q$
on the curve $\gamma$ is
$$\delta(p,q) := \frac{d(p,q)}{|p-q|}\ge1.$$
This is a Lipschitz constant for the inverse of the map
$\gamma$ (which itself has Lipschitz constant $1$).
The distortion of $\gamma$ is the maximum distortion between any points:
$$\delta(\gamma) := \sup_{p,q} \delta(p,q),$$
where the supremum is over $(\gamma\cross\gamma)\setm\Delta$, the set of
all pairs of distinct points on $\gamma$.  We can also consider a
restricted distortion, only considering opposite pairs:
$$\delta_\circ(\gamma) := \sup_{p\in\gamma} \delta(p,p^*).$$
Clearly $\delta_\circ(\gamma)\le\delta(\gamma)$ since the supremum
is over a subset; also, the supremum in $\delta_\circ$ is achieved,
though that for $\delta$ may not be.
\end{defn}

\begin{defn}
Given any real number $b\ge1$, the {\em distortion thickness}
of a curve $\gamma\subset\R^n$ is
$$\tau_b(\gamma) := \inf_{\delta(p,q)>b}|p-q|,$$
the minimum distance between points of distortion
more than $b$.
\end{defn}
The idea here is that on a smooth curve,
nearby points have distortion hardly greater than one, so by taking $b>1$
we eliminate these nearby pairs before taking the minimum distance.

We will define the remaining notions only for $C^{1,1}$ curves, which
have a well-defined, Lipschitz continuous tangent vector field,
and a curvature $\kappa$ almost everywhere.  A pair of points $(p,q)$
on the curve is called {\em symmetric} if the unit tangent vectors
$T_p$ and $T_q$ make equal angles with the chord $p-q$, that is,
$T_p\cdot(p-q)=T_q\cdot(q-p)$.  For instance, any two points on a
round circle are symmetric.  An ordered pair $(p,q)$ of distinct
points is said to be {\em self-critical} \cite{Lith1} if $q$ is a
critical point for the distance from $p$, that is $T_q\cdot(p-q)=0$.
A pair which is symmetric and self-critical is called {\em doubly
self-critical}.

For any $C^{1,1}$ curve $\gamma$, define
$c_1(\gamma)$ to be the minimum straight-line distance $|p-q|$
among self-critical pairs $(p,q)$,
and $c_2(\gamma)$ to be the minimum among doubly self-critical pairs.
Also, let $r(\gamma)$ be the minimum radius of curvature along the curve,
$r(\gamma):=1/\sup_\gamma{\kappa}$; this is positive because a $C^{1,1}$
curve has bounded curvature.

\begin{defn}
Given a real number $k>0$, the {\em curvature thickness} of a curve $\gamma$ is
$$\sigma_k(\gamma) := \min\bigl(kr(\gamma),c_2(\gamma)\bigr).$$
For $k=2$, this is (twice) the notion of thickness studied in \cite{Lith1}.
\end{defn}

\begin{note}
If $\gamma$ is not $C^{1,1}$, its curvature in some sense is unbounded,
so we set $\sigma_k(\gamma)=r(\gamma)=0$.
\end{note}

Since we are interested in minimizing the length of curves with 
(any one of our notions of) thickness prescribed, or dually, in maximizing
thickness for curves of prescribed length, it useful to have a scale-invariant
measure for this optimization problem.
\begin{defn}
The {\em ropelength} of a curve is its length divided by its
thickness.  This is the length of a similar curve scaled to have
thickness $1$ (so that it could be made out of a rope of diameter $1$).
Given a topological knot type $K$, the {\em ropelength of} $K$ is the
minimum ropelength of any curve $\gamma$ of type $K$.
\end{defn}
Intuitively, the ropelength of $K$ is the shortest length of rope
(of diameter one) which could be used to tie that knot.
Of course, there are several notions of ropelength corresponding
to our different notions of thickness, but we conjecture they all
agree (for reasonable values of $k$ and $b$), because a curve
of shortest ropelength should have nice structure.

\begin{note}
Most of the notions we have defined make sense also for
any union of immersed rectifiable curves, which need not be closed curves.
For a link of several components, we take $d(p,q)=\infty$
if $p$ and $q$ are on different components, so that such pairs
are always contenders in the definition of distortion thickness.
If the curves are not embedded, there are points
$p$, $q$ such that $d(p,q)>0$ while $|p-q|=0$.  Such a pair has infinite
distortion, and is doubly self-critical, so all our notions agree
that such a nonembedded curve has thickness zero.
\end{note}

\section{Distortion of Curves}

In order to understand the definition of distortion thickness $\tau_b$,
we need to explore the distortion of curves.  In particular,
we will see that every closed curve has distortion at least $\pi/2$,
so that $b=\pi/2$ is a natural choice.

Recall that the distortion of $\gamma$ is defined as a supremum
of $\delta(p,q) = d(p,q)/|p-q| \ge 1$.
We noted before that if $\gamma$ is not rectifiable or embedded
then we will have some $\delta(p,q)=\infty$.
The finiteness of $\delta(\gamma)$ is sometimes called a
{\em uniform chord-arc condition} in harmonic analysis \cite{Toro}.
We now note that this supremum may be infinite even for an embedded
rectifiable curve.

\begin{exmp}
Consider first a wedge $\gamma$ of external angle $\theta$.
If $p$ and $q$ are at equal distance $s$ from the corner, then
$|p-q|=2s\cos\frac\theta2$ so $\delta(\gamma)=\delta(p,q)=\sec\frac\theta2$.
Now take a plane curve $\gamma$, such as $x^2=y^3$, with a sharp cusp
at the origin.  Then as $p$ and $q$ approach the origin along $\gamma$
at equal distance $s$ from (and on opposing sides of) the cusp, we may
approximate a neighborhood of the cusp by narrower and narrower wedges,
and conclude that $\lim\delta(p,q)=\infty$.
\end{exmp}

However, for a $C^1$ curve, the distortion $\delta(p,q)$ approaches
$1$ whenever $p$ and $q$ approach a common point, so we can view $\delta$
as a continuous function on the compact space $\gamma\cross\gamma$,
and it achieves its (finite) supremum somewhere (away from the diagonal).

\begin{note}
O'Hara \cite{Ohara1} defined a family of energy functionals $E^p_j$ for
knots as $L^p$ norms
$$j E^p_j(\gamma) := \bigl\| |x-y|^{-j} - d(x,y)^{-j} \bigr\|_p$$
and observed that the distortion is a limiting value
$\delta(\gamma)=\exp(E^\infty_0(\gamma))$.
\end{note}

Gromov mentions in \cite{GromovFRM} that any closed rectifiable curve has
distortion at least $\pi/2$.  In fact, in \cite{GLP} he proves a
somewhat stronger result which includes the following proposition as a
consequence.  For completeness, we include a proof.

\begin{prop}
Given any closed rectifiable curve $\gamma\subset\R^n$,
the distortion satisfies $\delta(\gamma)\ge\delta_\circ(\gamma)\ge\pi/2$,
with equality if and only if $\gamma$ is a round circle.
\end{prop}

\begin{proof}  
Rescale so that $L(\gamma)=2\pi$, and parameterize the curve by arc-length
$s \in \R/2\pi\Z$.
Opposite points $p$ and $p^*$ have
arclength $d(p,p^*)=\pi$, so we want to prove for some $p$ that
$|p-p^*|\le2$.  For then $\delta(p,p^*)\ge\pi/2$ and we are done.

Consider the new curve $f(s) := p-p^*=\gamma(s)-\gamma(s+\pi)$ in $\R^n$.
Observe that $f$ is Lipschitz with speed at most two
(since by the triangle inequality
$|f'(s)|\le|\gamma'(s)|+|\gamma'(s+\pi)|= 2$ almost everywhere),
and that $f(s+\pi) = -f(s)$ for all $s$.

We want to show that $|f|\le2$ somewhere; suppose not, 
so that $f$ lies {\it outside} the closed ball of radius two in $\R^n$.
Then any arc of $f$ from $s$ to $s+\pi$
is an arc between antipodal points in $\R^n$ which avoids this closed
ball, so its length exceeds the distance
between antipodal points on a sphere of radius two, which is $2\pi$.
But since the parameterization of $f$ has speed at most two, the length
of this arc is at most $2\pi$, a contradiction.

Note that the only way to get equality $\delta=\pi/2$ is to have $|f|\ident2$,
with $f$ tracing out a great circle on this sphere.
Then $p$ and $p^*$ must always move in opposite directions, and
$p$ also traces out a round circle.
\end{proof}

Even if we had only proved this theorem for smooth curves, the
following lemma would extend the result to arbitrary curves;
it will also be useful in our discussion of distortion thickness.

\begin{lem}
In the $C^0$ topology on the space of curves, the distortion,
considered as a map to $[1,\infty]$, is lower semicontinuous.
\end{lem}

\begin{proof}
In the $C^0$ topology on parameterized curves (not necessarily
with constant speed parameterizations), the position of any point
$p=\gamma(t)$ is by definition a continuous function of $\gamma$.
Thus, for any points $p$ and $q$ the function $|p-q|$ is continuous in
the $C^0$ topology, and the arclength $d(p,q)$ along $\gamma$ is lower
semicontinuous.  Hence $\delta(p,q)$ is lower semicontinuous for each
pair of points.

Remember that a real-valued function $f$ on any topological space is
lower semicontinuous if and only if $f^{-1}\{x:x>y\}$ is open for all
$y\in\R$.  Thus it is clear that if every member of some family $\{f_\alpha\}$
of functions is lower semicontinuous, then $f:=\sup_\alpha f_\alpha$
is also lower semicontinuous.  But the distortion $\delta(\gamma)$
is the supremum of $\delta(p,q)$, a family of lower semicontinuous
functions, hence is itself lower semicontinuous.
\end{proof}

\begin{note}
This lemma, together with the fact that distortion blows up for
nonembedded curves, leads us to expect that a curve minimizing the distortion
should exist within each knot class.  Gromov \cite{GromovHED} and
O'Hara \cite{Ohara1} have independently observed that there are infinitely
many knot types with distortion less than some constant $M<100$.  Gromov asked
\cite[p.~114]{GromovFRM} if perhaps all knots have distortion under $100$.
In any case, distortion does not seem to be as useful for distinguishing
knots as O'Hara's original energy $E^1_2$ \cite{Ohara,FHW,KS-knot} has
been, or as thickness promises to be.
\end{note}

\begin{note}
On a sufficiently smooth curve $\gamma$, we can make some observations about
any pair $(p,q)$ realizing the maximum distortion.  First if we consider
the variation of $\delta(p,q)$ from a perturbation fixing $d(p,q)$, we
find that $(p,q)$ must be a symmetric pair, meaning the chord between them
makes equal angles $\theta$ with the tangents to $\gamma$ at either end.
If $p$ and $q$ are not opposite points on $\gamma$, then we can also
consider variations increasing or decreasing $d(p,q)$; this shows that
$\delta=\sec\theta$.  But if $q=p^*$ is opposite $p$, all we can conclude is that
$\delta\le|\sec\theta|$, which bounds $\theta$ near $\pi/2$.
Note that a doubly self-critical pair realizing the minimum self-distance
$c_2(\gamma)$ cannot also realize the maximum distortion, unless the points
are antipodal; some slightly more distant pair will have greater distortion.
\end{note}

\section{Distortion and a Conjecture in Integral Geometry}
Another approach to proving the proposition from the last section
would be to use Crofton's formula
and related notions from integral geometry, which relate the average
sizes of projected images of a curve to the length of the curve.

Recall that if $\gamma$ is a plane curve of length $L$,
and we consider the ($S^1$-worth of) projected images
of $\gamma$ to the line,
Crofton's formula says that the average length (with multiplicity)
of these projections is $\frac2\pi L$.
Since the multiplicity of the projected closed curve is at least 2 everywhere,
we find that the average diameter of the projections is at most $L/\pi$.
Thus, a plane curve of length $2\pi$ has width at most $2$ in some direction.

Our proposition follows immediately for plane curves.   If $\gamma$ has length
$2\pi$, orient it so that its width horizontally is at most $2$.
Then find a pair of opposite points $p$ and $p^*$ at the same height
(noting that the height difference changes sign).  Clearly $|p-p^*|\le2$,
while $d(p,p^*)=\pi$, so $\delta_\circ(\gamma)\ge\pi/2$.

Janse van Rensburg has suggested \cite{Buks} that any space curve
could be ``unfolded'' to a convex plane curve of the same length,
while never decreasing the chord distances $|p-q|$.  For instance,
if a curve touches a supporting hyperplane in two disjoint places,
we could reflect one intervening segment to lie below this plane.
This procedure certainly never increases the distortion.  If we
could prove that it converged to a planar curve
(perhaps considering only the polygonal case)
then we would have a new proof of our proposition, except that it would
seem difficult to show circles are the only curves of distortion
$\delta=\pi/2$.

The argument sketched above for plane curves could be applied directly to
curves in higher dimensions if the following conjecture is true.
Again we would find a pair of opposite points on the original curve
that have the same height in the direction of the projection,
and they would then have distortion at least $\pi/2$.

\begin{conj}
If $\gamma$ is a curve in $\R^n$ of length $L$, then there is
some orthogonal projection to $\R^{n-1}$ in which the image
of $\gamma$ has diameter at most $L/\pi$.
\end{conj}

\begin{note}
For space curves in $\R^3$ of length $L$,
there are two standard analogues of Crofton's formula.
The first deals with projections to a line, and says that
the average length (with multiplicities), over the sphere of projections,
is $L/2$.  Thus some projection to a line has diameter at most $L/4$.
The second deals with projections to a plane, and says that
the average length of these, over the sphere of possible projections,
is $\frac\pi4 L$.  Thus some projection to a plane has diameter
at most $\frac\pi8 L$, but this is not as good as our conjectured
$L/\pi$.
\end{note}

\section{The Distortion Thickness}

Recall that the distortion thickness of a curve
measures the least self-distance among pairs of points with large distortion:
$$\tau_b(\gamma) := \inf_{\delta(p,q)>b}|p-q|.$$
Clearly, as we increase $b$, $\tau_b(\gamma)$ is nondecreasing,
since we take the infimum over smaller sets.
For $b\le1$, all pairs $(p,q)$ are in contention,
except those connected by a straight arc, so $\tau_1=0$ for all closed curves.
For $b\ge\frac\pi2$, the infimum may be over the empty set (as it
would be for a circle), so the thickness may be infinite.
For $b=\frac\pi2$, although a circle has infinite thickness,
no other closed curve does.
Thus we concentrate on $b\in(1,\frac\pi2]$.

\begin{note}
If $b<\delta(\gamma)$, then the distortion ropelength $L/\tau_b$ is
always at least $2\delta(\gamma)\ge\pi$.  To check this, take $p$ and $q$
with $\delta(p,q)$ close to $\delta(\gamma)$.  Of course $d(p,q)\le L/2$
and $|p-q|\ge \tau_b$, so dividing gives us $\delta(p,q)\le L/2\tau_b$.
\end{note}

\begin{prop}
The distortion thickness $\tau_b$ is upper semicontinuous
in the $C^0$ topology on the
set of closed rectifiable curves, for any fixed $b$.
\end{prop}
\begin{proof}
The claim is that if a sequence of curves $\gamma^k$ approaches a
limit $\gamma^0$,
then their thicknesses $\tau^k:=\tau_b(\gamma^k)$ satisfy
$$\limsup\tau^k \le \tau^0.$$
If not, write $\limsup\tau^k = \tau^0+3\epsilon$, for some
$\epsilon > 0$; passing to a subsequence, we may
assume that $\tau^k>\tau^0+2\epsilon$ for all $k$.

Now, the infimum in the definition of $\tau^0$ may not be realized,
but we can certainly find a pair of points $(p,q)$ which come close to
realizing $\tau^0$.  That is, $|p-q|_0$, their straight-line distance
on $\gamma^0$, is within $\epsilon$ of $\tau^0$, while $\delta^0(p,q)$,
their distortion on $\gamma^0$, is greater than $b$.  So we have
$$\tau^k>|p-q|_0+\epsilon, \qquad \delta^0(p,q)>b.$$
Since $|p-q|_k\approaches|p-q|_0$, but $\tau^k$ stays greater than
this distance, we must have that $(p,q)$ is not in contention for
$\tau^k$ for all large enough $k$, that is, $\delta_k(p,q)\le b$.
On the other hand, the semicontinuity of $\delta$ means that
$\liminf\delta_k(p,q)\ge\delta_0(p,q)>b$.
This contradiction completes the proof.
\end{proof}

\begin{note}
The length of curves is lower semicontinuous (but not continuous)
in this $C^0$ topology.  Although we have discussed the existence of an
arclength parameterization for each curve, we do not assume when discussing
limits that constant-speed parameterizations are used.  For instance,
consider curves obtained by replacing one side of a square by finer and finer
zigzags of twice the length.  If each is parameterized by arclength,
they approach a limit, which traces out the square at varying speed.
\end{note}

\begin{cor}
In the $C^0$ topology, the distortion ropelength of curves is
lower semicontinuous.
\end{cor}
\begin{proof}
It is the quotient of a lower semicontinuous function (arclength) by an
upper semicontinuous function (thickness).
\end{proof}

Let us consider some example curves.  First take $\gamma$ to be the unit circle.
Given any $\phi\in[0,\frac\pi2)$, any pair of points at angle $2\phi$ is
at distance $|p-q|=2\sin\phi$, so it has distortion $\phi/\sin\phi$.
Thus if $b=\phi/\sin\phi\in[1,\frac\pi2)$ the thickness $\tau_b(\gamma)$
will be $2\sin\phi$.  (The infimum here is not realized.)

Next, recall that the distortion of a wedge of (external) angle $\theta$ is
$\sec\frac\theta2$.
Given $b$, set $\theta=2\sec^{-1}b$; if a curve $\gamma$ has any corner
sharper than this, its thickness $\tau_b(\gamma)$ will be zero.
However, on a polygonal curve with no corners sharper than this,
no pair of points on adjacent segments will be in contention for $\tau_b$.

Suppose $\gamma$ bounds the convex hull of the unit circle and an exterior
point (placed so that the angle at this point is $\theta$).  This curve has
thickness $\tau_b=0$ if $b<\sec\frac\theta2$, but if
$b\ge\sec\frac\theta2$ its thickness is the same as that of the circle.
This shows that the thickness is not $C^0$-continuous.

\begin{prop} \label{tauctwo}
Suppose $\gamma$ is a $C^{1,1}$ curve.  If the infimum defining $\tau_b$
is attained, this is at a doubly self-critical pair, so $\tau_b\ge c_2$.
If not, but $\tau_b<\infty$, there is a pair
of points with $\delta(p,q)=b$ and $|p-q|=\tau_b$.
\end{prop}
\begin{proof}
If the infimum is attained at $(p,q)$ with $\delta(p,q)>b$,
then all nearby pairs on $\gamma$ also have $\delta>b$
and are in contention for $\tau_b$,
so since $|p-q|$ is a minimum, $(p,q)$ must be a doubly self-critical pair.
Otherwise, take a minimizing sequence $(p_n,q_n)$, and pass to a
convergent subsequence with some limit $(p,q)$.
Because the sequence was minimizing, $|p-q|=\tau_b(\gamma)$, with
$\delta(p,q)\ge b$.  But if $\delta(p,q)>b$, this pair would
achieve the minimum, so instead $\delta(p,q)=b$.  (Note also that this
implies $d(p,q)=b\tau_b$.)
\end{proof}

\begin{note}
The minimum doubly self-critical distance $c_2$ can be achieved
at a pair of arbitrarily small distortion $\delta>1$, as if we put
small hooks at the ends of a straight segment before continuing
with a huge loop.  Thus it is not clear in general that $\tau_b\le c_2$.
Perhaps, however, the arc between the points in a self-critical pair
always has large distortion somewhere.  Hence we might consider a refined
definition of $\tau_b$ which allows in the infimum any $(p,q)$ for which
a subarc has distortion $\delta>b$.
\end{note}

\section{The Curvature Thickness}
Recall that
$\sigma_k(\gamma) := \min\bigl(kr(\gamma),c_2(\gamma)\bigr),$
which is the smaller of the minimum doubly self-critical distance $c_2$
and the minimum radius of curvature $r$ scaled by the factor $k$.
This definition makes sense for any $k\ge0$, but we are most interested
in $k\in[1,2]$.

The special case of this definition with $k=2$ was introduced in \cite{Lith1}.
This thickness $\sigma_2$ is characterized there as the maximum diameter of
a tube around a $C^2$ curve $\gamma$ for which the normal exponential map
is an embedding.  We will discuss this and similar notions in
Section~\ref{tubes}.

The further main results of \cite{Lith1} deal also with $\sigma_2$
for $C^2$ curves and, translated into our notation, say that
\begin{enumerate}
\item If $\gamma$ has ropelength $L/\sigma_2<\frac{n\pi}2$ for some
integer $n$, then
there is some $n$-gon isotopic to $\gamma$ (in the same knot class).
\item Thus, any nontrivial knot has ropelength at least $\frac{5\pi}2$.
\item If $\gamma$ is a $C^2$ curve then its ropelength is at least twice its
distortion: $L(\gamma)/\sigma_2(\gamma) \ge 2\delta(\gamma)$.
\item If $\gamma$ is a $C^2$ curve, then $c_1(\gamma)\ge\sigma_2(\gamma)$.
In other words, $c_1(\gamma)=c_2(\gamma)$ when either is less than $2r(\gamma)$
\end{enumerate}

\begin{note}
This last result shows that $\sigma_k=\min(kr,c_1)$ for any $k\le2$,
so that we could have defined $\sigma_k$ in terms of $c_1$ instead of $c_2$.
\end{note}

\begin{note}
O'Hara \cite{Ohara1} attributes the definition of $c_1(\gamma)$
to Kuiper, and shows that if $\gamma$ has $L/c_1<n$ then it is isotopic to
an $n$-gon.  This is a weaker conclusion than the similar
result of \cite{Lith1}, but requires no curvature bound.
\end{note}

\begin{note}
One of the authors has pointed out \cite{Sul-ima}
that the discretization for thickness suggested by Stasiak \cite{Sta}
(and presumably used in numerical work reported in \cite{Katritch})
actually discretizes
$\sigma_1$.  Here we approximate a smooth curve $\gamma$ by a polygon
with short edges of nearly equal length.  Now look at solid cylinders
of equal diameter around each edge.  Increase the diameter until some
nonadjacent cylinders first intersect; this gives an approximation
of $\sigma_1(\gamma)$.
\end{note}

It is clear that $\sigma_k=0$ for any curve with a sharp corner.
Thus it is not possible to have any useful bound of the form
$\sigma_k\ge C \tau_b$.
However, we can bound $\tau_b$ below in terms of $\sigma_k$,
because a curvature bound implies that short segments of a curve
have small distortion.  We use a lemma,
which improves one from \cite[p.~242]{Ohara} by a factor of two,
and is in fact just a special case of
Schur's Theorem (see \cite[p.~36]{Chern}).
\begin{lem}
If $\kappa(\gamma)\le 1$, then for any points $p$ and $q$
on $\gamma$ with $d(p,q)=d\le2\pi$, we have $|p-q|\ge2\sin d/2$, or
equivalently,
$$\delta(p,q)\le\frac {d/2}{\sin d/2}.$$
\end{lem}
\begin{proof}
Let $x$ be the point halfway along the arc from $p$ to $q$.
Orient the curve so that the tangent vector is vertical at $x$.
At arclength $s$ from $x$, the curve
has turned less than angle $s$, so the vertical component of the
unit tangent vector is still at least $\cos s$.  Thus out to arclength $d/2$,
the vertical component of the difference vector $x-p$ or $q-x$ is at least
$\int_0^{d/2}\cos s\,ds=\sin d/2$.  Thus $|p-q|\ge2\sin d/2$, as desired.
\end{proof}

\begin{defn}
For any $b\in[1,\frac\pi2]$, define $k_b$ to be the unique solution
of the equivalent equations
$$b = \frac{\arcsin{k_b/2}}{k_b/2},\qquad \frac{bk_b/2}{\sin bk_b/2}=b.$$
Note that $k_b$ increases monotonically from $0$ to $2$ as $b$ increases
from $1$ to $\frac\pi2$, and thus $bk_b\le\pi$.
\end{defn}

\begin{thm} \label{distocurv}
For any $1\le b\le\frac\pi2$ and any closed curve $\gamma$, we have
$$\tau_b(\gamma) \ge \sigma_{k_b}(\gamma).$$
\end{thm}
\begin{proof}
Write $k=k_b$.
If $\sigma_k=0$ (as when $\gamma$ is not $C^{1,1}$) there is nothing to prove.
Otherwise, rescale so that $\sigma_k(\gamma)=k$.  Then
we have $c_2(\gamma)\ge k$ and $\kappa\le 1$ everywhere on $\gamma$.
Now apply Proposition~\ref{tauctwo};
if the infimum defining $\tau_b$ is attained,
$\tau_b\ge c_2\ge k$, and we are done.

Otherwise, the proposition gives us points $p$ and $q$ with
$|p-q|=\tau_b$ and $d(p,q)=b\tau_b$.
If $\tau_b<k$ we want to use the lemma to derive a contradiction.
Certainly then $d(p,q)<kb\le\pi$ so the lemma applies, giving
$$b=\delta(p,q)\le\frac{b\tau_b/2}{\sin b\tau_b/2}.$$
But $\frac x{\sin x}$ is an increasing function, so since $b\tau_b<bk$, we
get $b<\frac{bk/2}{\sin bk/2}$ contradicting the definition of $k_b$.
\end{proof}

In particular, we have $\tau_{\frac\pi3}\ge\sigma_1$ and
$\tau_{\frac\pi2}\ge\sigma_2$.
Since $\tau_b(\gamma)$ and $\sigma_k(\gamma)$ are nondecreasing in
$b$ and $k$ respectively, it follows that $\tau_b\ge\sigma_k$
whenever $b>2\arcsin(k/2)/k$.  We can also get bounds for other $b$ and $k$.

\begin{cor}
For any $b\in[1,\frac\pi2]$, $k\ge k_b$ and any curve $\gamma$,
$$\tau_b(\gamma) \ge \frac{k_b}k \sigma_k(\gamma).$$
In particular, $\tau_\frac\pi3 \ge \sigma_k/k$ for $k\ge1$.
\end{cor}
\begin{proof}
Apply the theorem to $b$ and $k_b$, and use the fact that by definition of
$\sigma_k$, we have
$\sigma_{k_b}\le\sigma_k\le\frac{k}{k_b}\sigma_{k_b}$.
\end{proof}

\begin{note}
Combining the result $\tau_\frac\pi2 \ge \sigma_2$ with our
previous note about distortion ropelength, we recover the
result from \cite{Lith1} that $L/\sigma_2\ge2\delta(\gamma)$.
In fact, their proof is similar in spirit to our proof
of Theorem~\ref{distocurv}.
\end{note}

\begin{note}
If $\gamma$ is the unit circle, then for any $b$ we have
$\tau_b(\gamma)=\sigma_{k_b}(\gamma)=k_b$, so the theorem
is sharp.
\end{note}

\section{Thickness Notions Related to Normal Tubes}\label{tubes}

Other authors have suggested various notions of thickness defined
as the largest diameter of a normal tube before some property
fails.  From \cite{Lith1}, we have seen that $\sigma_2$ is the infimum of
diameters for which the normal exponential map fails to be an embedding.

Diao {\sl et al.}~\cite{DEJvR} suggest that we require the
$\frac{d}2$-neighborhood to be a solid torus and then look at the
intersection of this neighborhood with each normal plane (assuming $\gamma$
is $C^1$).  If further we ask that the zero-component of each such intersection
be a meridian disk in the solid torus, intersecting $\gamma$ only at the
original point, then the first $d$ for which
this fails is their $C^1$-continuous measure of thickness,
which we call $t_D(\gamma)$.

O'Hara \cite[p.~60]{Ohara2} suggests
defining thickness simply as the infimum of all
$d$ such that the $\frac{d}2$-neighborhood of $\gamma$ is not topologically
a solid torus; we will call this $t_O(\gamma)$.
However, the notion he actually uses is different.  Let $t_B$ be
the infimal diameter of balls in space whose intersection with $\gamma$
is not a single arc (and thus is disconnected or all of $\gamma$).
O'Hara's implicit claim is that $t_B\le t_O$.\footnote{In fact, for any
diameter greater than $t_O$, he
claimed to get a ball with disconnected intersection with $\gamma$.
This can never be true if $\gamma$ is a circle, so we have tried to
correct for this in our definition of $t_B$.}
Certainly we can have $t_B<t_O$, as for a polygon, where $t_B=0$.

It is clear from these definitions that $\sigma_2 \le t_D \le t_O$,
and that $t_D \le c_1$.  Combining these with the final result of
\cite{Lith1} shows that $t_D=c_1=c_2$ whenever these are less than $2r(\gamma)$.

For the ellipse $(x/a)^2+(y/b)^2=1$ with $a<b$, we have $c_2=2a=t_O$,
$r=a^2/b<a$, and $c_1=t_D$ with $\sigma_2=2r<t_D<t_O$.
(Explicit values of $\tau_b$ here would involve elliptic integrals.)
Suppose we now take the upper half of this ellipse, and attach it
in a $C^{1,1}$ fashion to a large loop in the lower halfplane;
this new curve has $c_2=2a<<t_O$.

Now consider a helix of radius $\rho$ and pitch $p$, the intersection
of a cylinder of radius $\rho$ with the helicoid $z=pr\cos\theta$.
(To get a closed curve, we could splice a long piece of this
helix into a huge loop.)  For steep pitch,
this has no doubly self-critical points, but if
$p^2\le-\min\bigl(\frac{\sin x}x\bigr)$, we have $c_2<3\rho$, approaching $0$
as the pitch decreases.  If the pitch is shallow enough, $c_2\le2\rho$ and
then $t_O=c_2$, but for an intermediate range of pitches
(when $c_2>2\rho$ but is still finite), we again find $t_O>>c_2$,
because by the time the normal tube is thick enough to see the 
self-critical points, it already extends across the axis of the cylinder,
so it is still homotopy equivalent to the helix.

\section{Shortest Curves in a Knot Class}
We have made some progress in understanding the structure of curves
minimizing the distortion ropelength, and believe the following:
\begin{conj}
For each nontrivial knot type $K$ there exists a shortest curve
$\gamma \subset \R^3$ of distortion thickness $1$ in that knot type;
that is, each knot type has a curve $\gamma$ of shortest distortion ropelength.
Any such $\gamma$ has bounded curvature, and thus is of
smoothness class $C^{1,1}$.  Away from doubly self-critical points realizing
the thickness, $\gamma$ must be straight.  If we scale to make the distortion
thickness $\tau_{\frac\pi2}(\gamma)=1$, then the curvature of $\gamma$ is bounded by $1$,
so that also the curvature thickness $\sigma_k(\gamma)=1$ for all $k\ge1$.
Finally, for any $b\ge\pi/3$, we have $\tau_b(\gamma)=1$, so all our
notions of thickness agree for a ropelength minimizer $\gamma$.
\end{conj}

\begin{note}
It seems that the ropelength minimizer for the Hopf link must be the obvious
candidate: two unit circles in perpendicular planes passing through each
others' centers.  For this link, $r=c_1=c_2=1$, so $\sigma_k=1$ for
any $k\ge1$.  Also, the distortion thickness for any $b\ge\frac\pi3$
is $\tau_b=1$.  If we extend the definitions of $t_D$, $t_O$ and $t_B$
to links in the obvious ways, these also all equal $1$.
Presumably, in the ropelength minimizer for the connected sum of
two Hopf links, the middle component (linking both others) will
be a stadium curve, built of two semicircles and two straight segments;
this shows we cannot expect minimizers to be $C^2$.
\end{note}


\begin{thebibliography}{[DEJvR]}

\bibitem[Che]{Chern}
S.~S. Chern.
\newblock Curves and surfaces in euclidean space.
\newblock In S.~S. Chern, editor, {\em Studies in Global Geometry and
  Analysis}, pages 16--56. Math. Assoc. Amer., 1967.

\bibitem[DEJvR]{DEJvR}
Yuanan Diao, Claus Ernst, and E.~J. Janse~van Rensburg.
\newblock Energies of knots.
\newblock Preprint, 1996.

\bibitem[FHW]{FHW}
Michael~H. Freedman, Zheng-Xu He, and Zhenghan Wang.
\newblock On the {M\"obius} energy of knots and unknots.
\newblock {\em Annals of Math.} {\bf 139}(1994), 1--50.

\bibitem[Gro1]{GromovHED}
Mikhael Gromov.
\newblock Homotopical effects of dilatation.
\newblock {\em J. Differential Geometry} {\bf 13}(1978), 303--310.

\bibitem[Gro2]{GromovFRM}
Mikhael Gromov.
\newblock Filling {R}iemannian manifolds.
\newblock {\em J. Differential Geometry} {\bf 18}(1983), 1--147.

\bibitem[GLP]{GLP}
Mikhael Gromov, J.~Lafontaine, and P.~Pansu.
\newblock {\em Structures m\'etriques pour les vari\'et\'es riemanniennes}.
\newblock Cedic/Fernand Nathan, Paris, 1981.

\bibitem[JvR]{Buks}
E.~J. Janse~van Rensburg.
\newblock Personal communication, June 1996.

\bibitem[KBM{\etalchar{+}}]{Katritch}
Vsevolod Katritch, Jan Bednar, Didier Michoud, Robert~G. Scharein, Jacques
  Dubochet, and Andrzej Stasiak.
\newblock Geometry and physics of knots.
\newblock {\em Nature} {\bf 384}(November 1996), 142--145.

\bibitem[KS]{KS-knot}
Robert~B. Kusner and John~M. Sullivan.
\newblock {M\"obius} energies for knots and links, surfaces and submanifolds.
\newblock In Willam~H. Kazez, editor, {\em Geometric Topology}, pages 570--604.
  Amer. Math. Soc./International Press, 1997.
\newblock Proceedings of the {G}eorgia {I}nt'l {T}opology {C}onference, August
  1993.

\bibitem[LSDR]{Lith1}
Richard~A. Litherland, Jon Simon, Oguz Durumeric, and Eric Rawdon.
\newblock Thickness of knots.
\newblock Preprint, 1996.

\bibitem[Oha]{Ohara}
Jun O'Hara.
\newblock Energy of a knot.
\newblock {\em Topology} {\bf 30}(1991), 241--247.

\bibitem[Oha1]{Ohara1}
Jun O'Hara.
\newblock Family of energy functionals of knots.
\newblock {\em Topology Appl.} {\bf 48}(1992), 147--161.

\bibitem[Oha2]{Ohara2}
Jun O'Hara.
\newblock Energy functionals of knots {II}.
\newblock {\em Topology Appl.} {\bf 56}(1994), 45--61.

\bibitem[Sta]{Sta}
Andrzej Stasiak.
\newblock Lecture at the AMS special session on {\em Physical Knot Theory},
  Iowa City, March 1996.

\bibitem[Sul]{Sul-ima}
John~M. Sullivan.
\newblock Lecture at the IMA workshop on {\em Topology and Geometry in Polymer
  Science}, June 1996.

\bibitem[Tor]{Toro}
Tatiana Toro.
\newblock Lecture at the IMA workshop on {\em Topics Related to Nonlinear PDE},
  March 1996.

\end{thebibliography}

\newcommand{\etalchar}[1]{$^{#1}$}

\end{document}